# A GEOMETRICAL LINK BETWEEN THE CIRCLE AND SEXAGESIMAL SYSTEM

By

Jaime Vladimir TORRES-HEREDIA JULCA
(University of Geneva, Switzerland)

## 1. - SUMMARY

This paper presents a simple geometrical fact which could relate to the history of mathematics and astronomy. This fact shows a natural link between the circle and the multiples of 6 and it makes it possible to obtain a simple representation of the 12 months of the year, the 24 hours of the day, the 30 days (average number) of the month and the 360 days (approximate number) of the year, which brings us closer to the sexagesimal division of time. Moreover this representation reminds one of the movement of the planets around a centre.
   Using this fact one will be able also to find geometrically the principal divisor of number 60, to represent numbers in base 60 with a kind of abacus or calculation table and to make a division of the circle into 6 and 12 equal parts. Afterwards one will be able to obtain a division in 360 unequal parts but relatively close to one another, and the goal isn't precisely to obtain an optimal division of the circle in 360 equal parts but to prove that the idea to divide the circle in 360 equal parts can subsequently be suggested by these geometrical facts that have been showed.
   In this article the author will not answer the following questions:
   a) What is the origin of the sexagesimal system?
   b) By which way could one manage to adopt the sexagesimal system starting from the knowledge of the facts exposed in this article and starting from the knowledge of the astronomical data?
   These questions could be treated, using information of this article, by the readers or later on by the author.

## 2. - INTRODUCTION

The sexagesimal division of time is widely used just like the sexagesimal division of the circle to measure the angles. However, one often wonders why one adopted these divisions. As a matter of fact, it is rather clear that the adoption of this system is linked up, among other things, to the need to represent astronomical data like the duration of the year (about 360 days).
   However, as this article will show, the sexagesimal division of the circle and time has close links with a simple geometrical fact that relates to the circle and that is independent of the astronomical data.
   Perhaps if this geometrical fact, with the need to represent the 365 days of the year, hasn't contributed to the adoption of the sexagesimal system, employed to divide the circle and to represent time, at least it can be used to show a natural geometrical link between the circle and the sexagesimal division.

## 3. - A NATURAL LINK BETWEEN THE CIRCLE AND NUMBERS 6, 12, 24, 30, 60 and 360

A natural link between the circle and numbers 6, 12, 24, 30, 60 and 360 appears clearly when one observes how a given disc is surrounded by other discs of the same radius successively, by building juxtaposed orbits around the initial disc.

This is the initial disc:

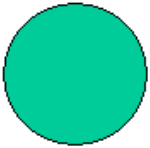

One surrounds this first disc with 6 discs of the same radius. The 6 new discs will compose an "orbit":

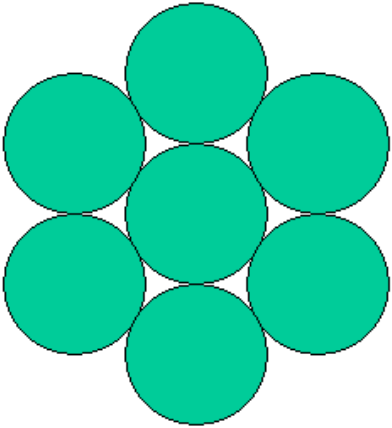

One repeats the process and a second orbit of 12 discs is added around to the first orbit:

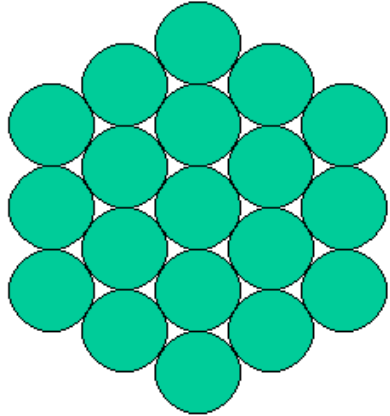

One can notice that when one adds orbits of discs of the same radius around the initial one, "equilateral triangles" made of discs build up around the centre. These triangles are distinguished below by different colours:

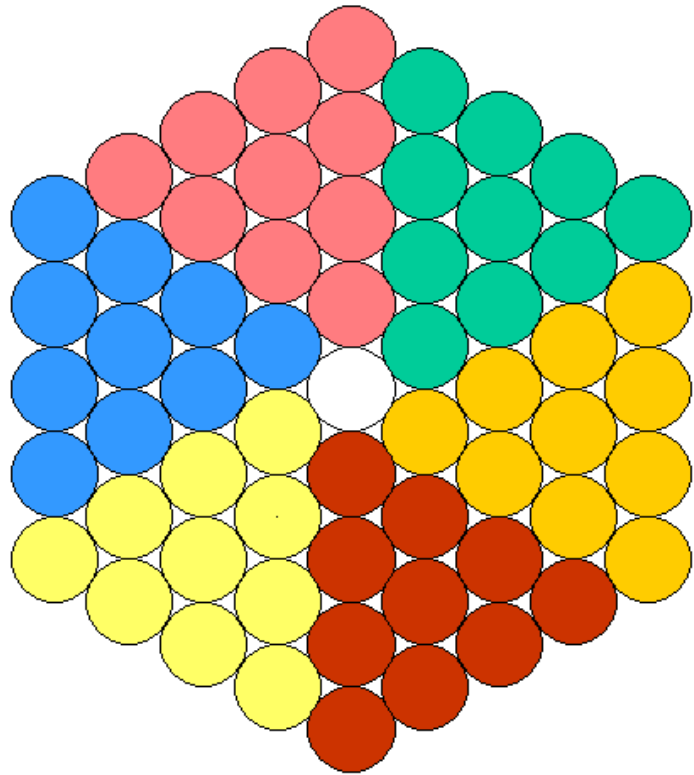

The more orbits there are, the bigger the 6 "equilateral triangles" made of discs get. One can also notice that one can find the number of discs of a given orbit according to the rank of the orbit. If one sets down that r indicates the rank of an orbit and that r = 1 for the first orbit,

i.e. for the first 6 discs around the initial disc, then number N of discs for a given orbit is:

$$N = 6 \times r$$

The number of discs of an orbit is therefore a multiple of 6. And if this formula is correct, the orbit of rank 60 will have 60x6 discs, i.e. 360.

This observation shows us that numbers 6, 12, 24, 30, 60 and 360 are already linked in a natural way to the disc and to the circle, because if one surrounds a given disc by the way suggested, one will obtain, at a certain point, 6, 12, 24, 30, 60 and 360 discs around the initial circle.

It is clear that there is an infinity of possible orbits since there is an infinity of multiples of 6. And one can choose the $60^{th}$ rank as representation of the days of the year because the advantage of number 360 is that it is also multiple of 5 and of 10.

Furthermore, if one counts the discs that are in the coloured triangles of the previous figure, one will have a total of 60, each "equilateral triangle" having 10 discs. Thus there are not only 60 discs but also the representation of 6x10 = 60. And the symmetries of this figure also make it possible to divide it into 2 and 3 equal parts. And as each of the 6 "equilateral triangles" can be divided into 2 and 5 equal parts (because they contain 10 discs), one can say that number 60 is divisible by 2, 3 and 5 and also by 4. In other words, this figure shows the principal divisors of number 60. And thanks to it one can also obtain another geometrical representation of number 360 with the total of the discs coloured below:

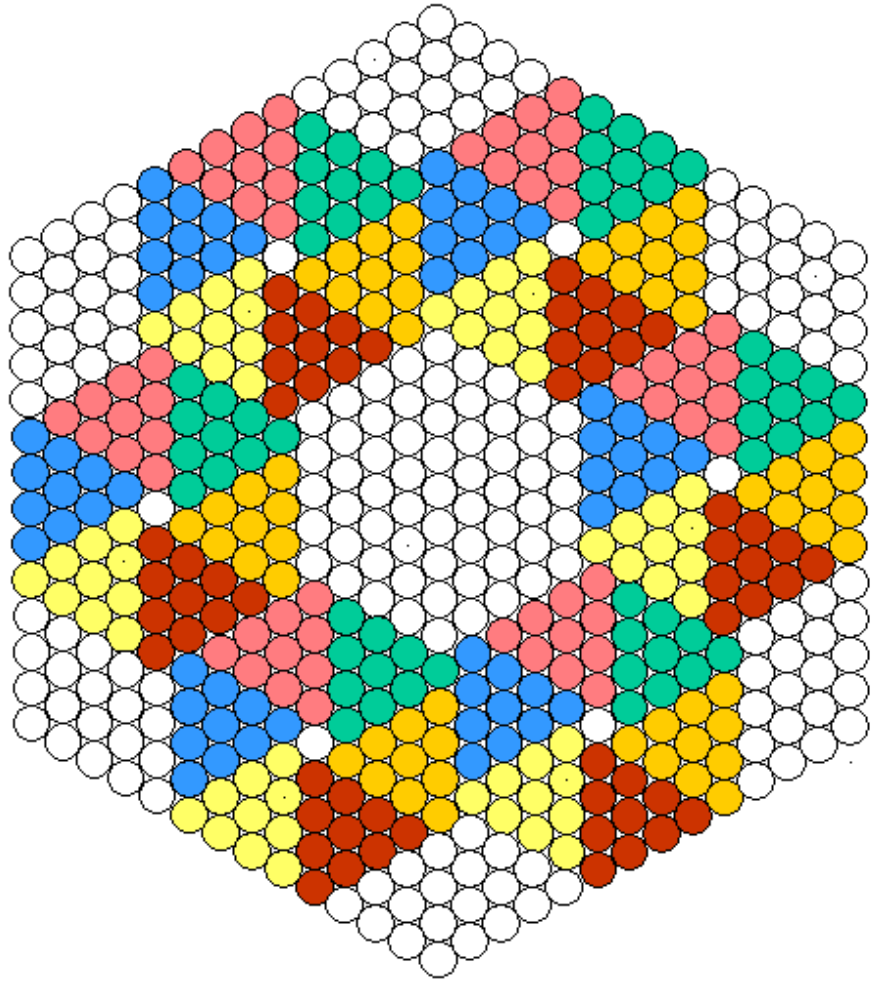

One can also represent number 360 like that, always with the total of the coloured discs:

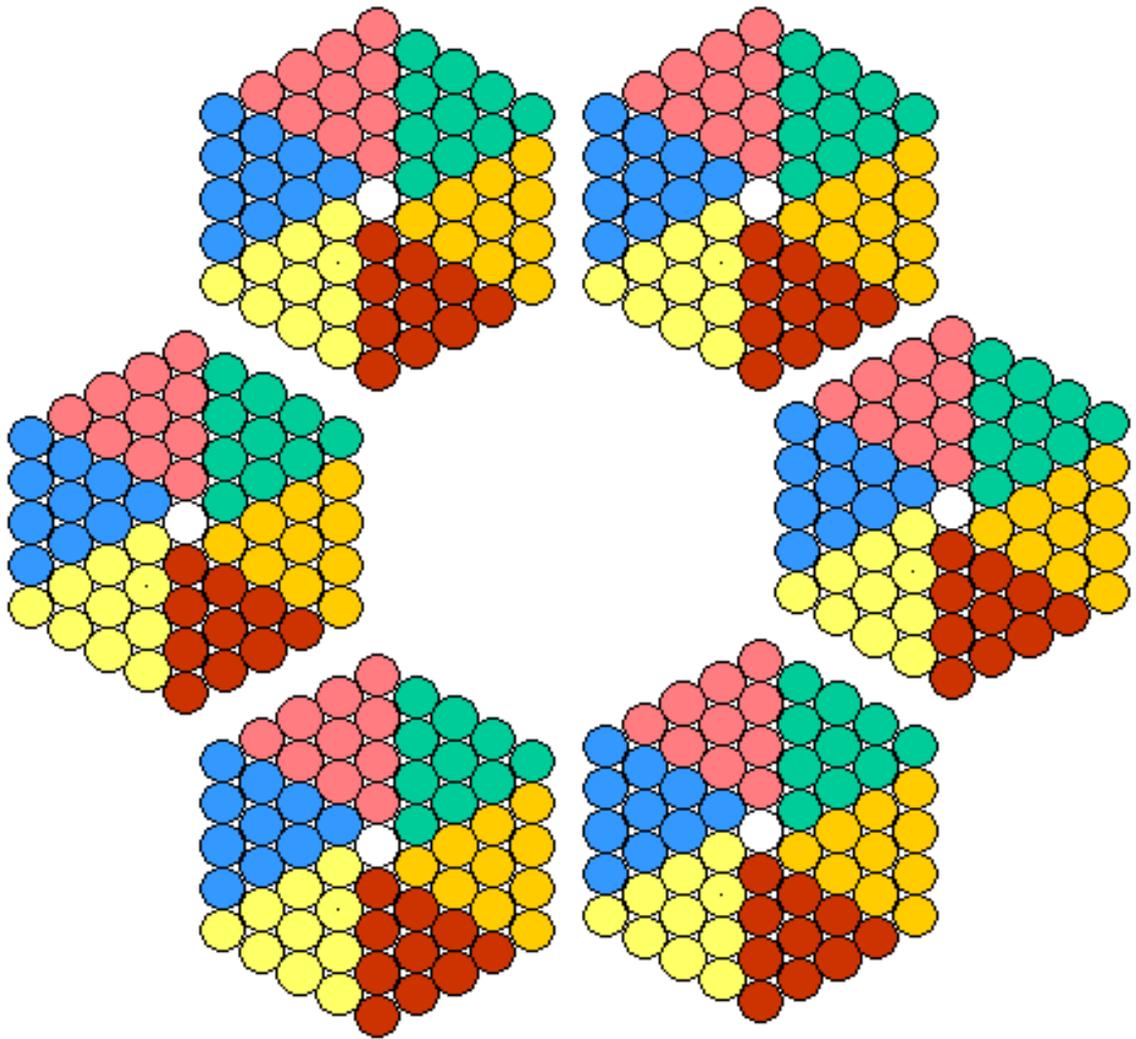

Moreover, if one chooses the sexagesimal system to measure time, one will be able to represent the number of months of one year with an orbit of 12 discs, the number of hours of one day with an orbit of 24 discs, and the average number of days of one month with an orbit of 30 discs, the whole thing in a way that reminds one of the movement of the planets around a centre.

Here is the representation of the twelve months of the year with brown discs:

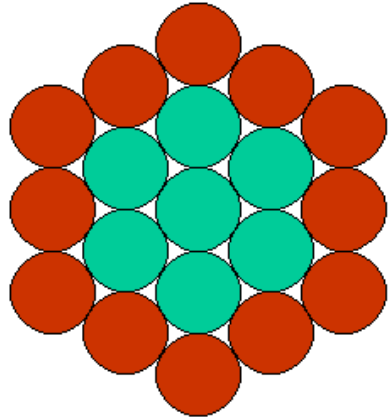

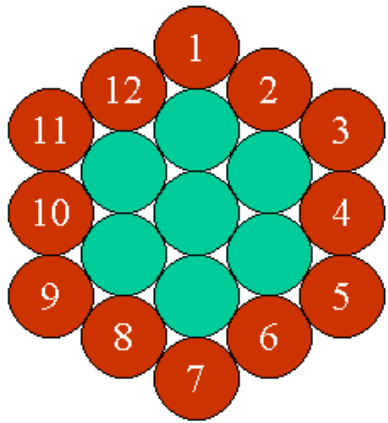

The 24 hours of the day with brown discs:

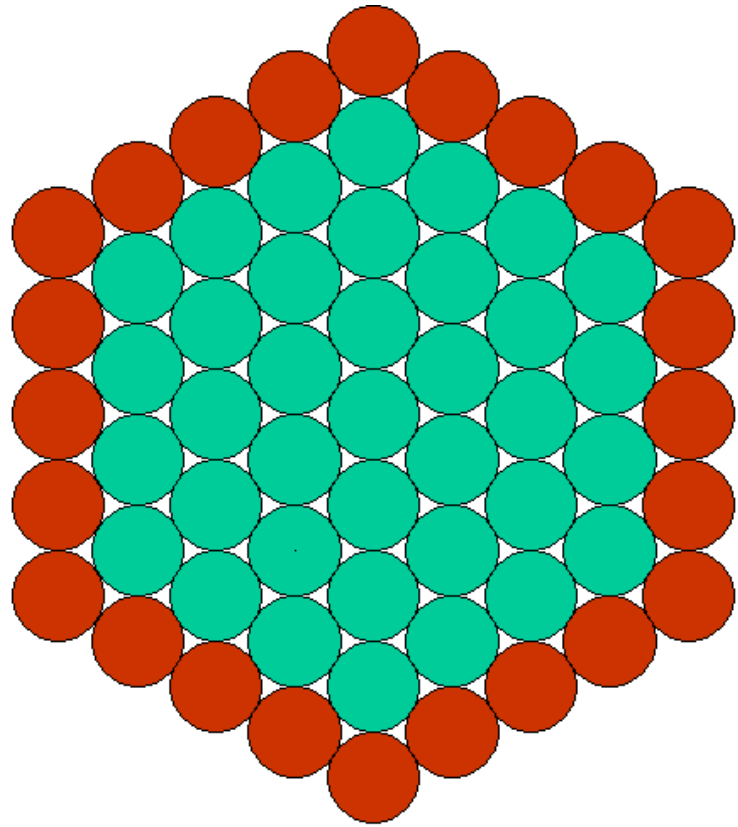

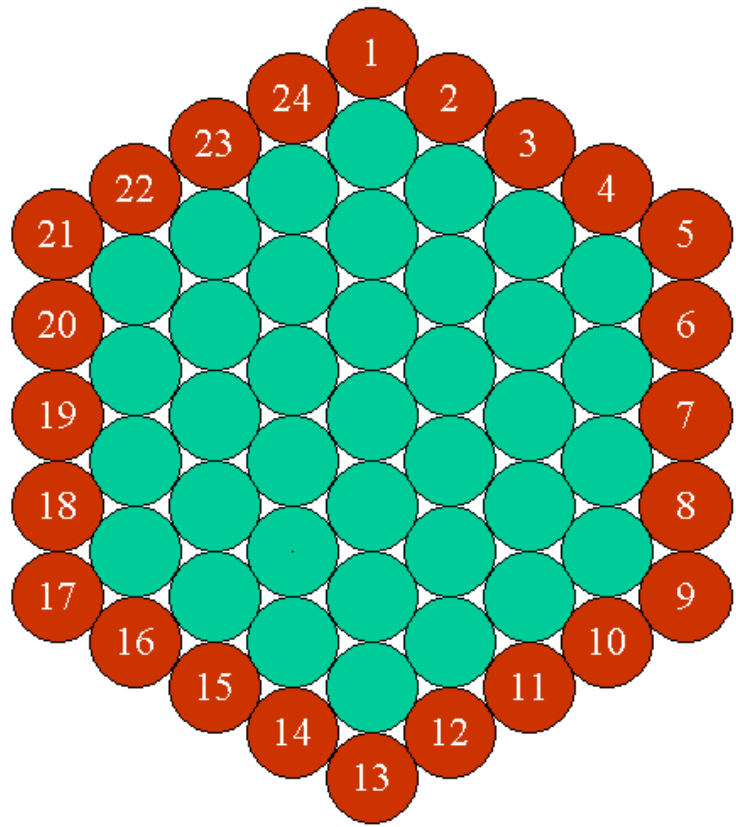

The 30 days of one month (average number):

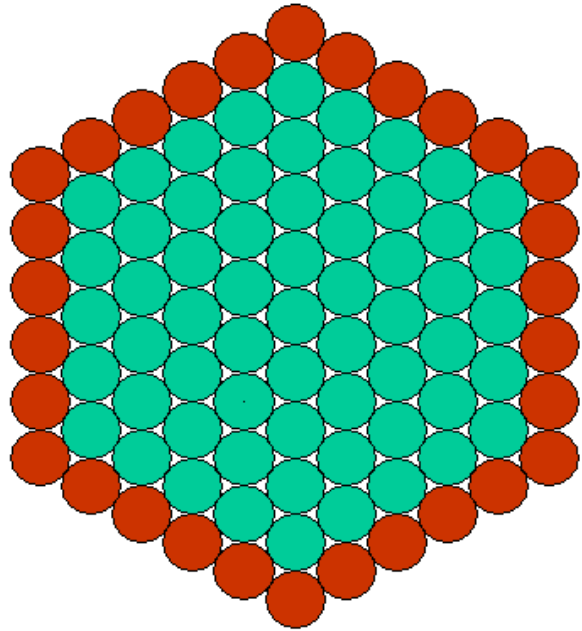

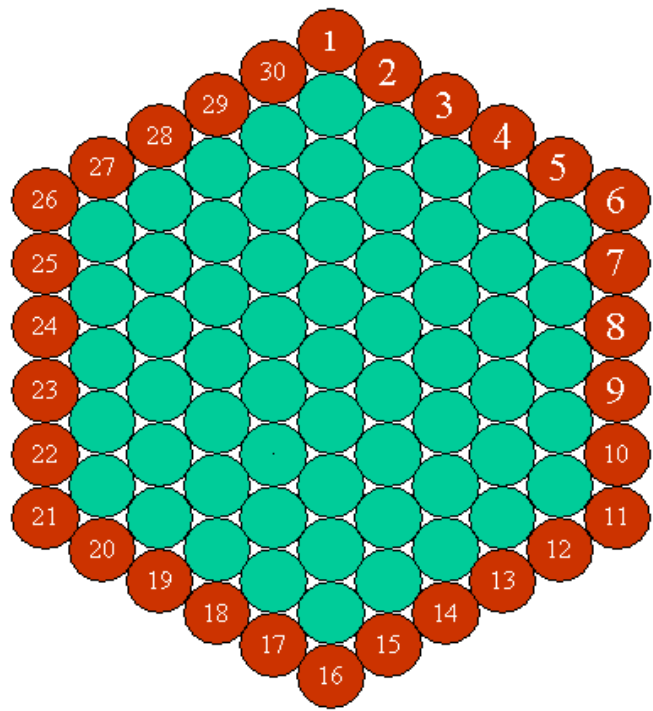

One could also represent by this way the approximate number of days of the year, i.e. 360.

This way of surrounding the circles enables us to foresee a division of the circle in 360 parts by the fact that number 360 is linked in a natural way, as we have seen, to the disc and the circle.

### 4. - ABACUS OR CALCULATION TABLE IN BASE 60

By the fact that one can represent number 60 as we have seen previously, one can build a kind of calculation table or abacus that will make it possible to represent the numbers in base 60.

Here is the principle of this table:

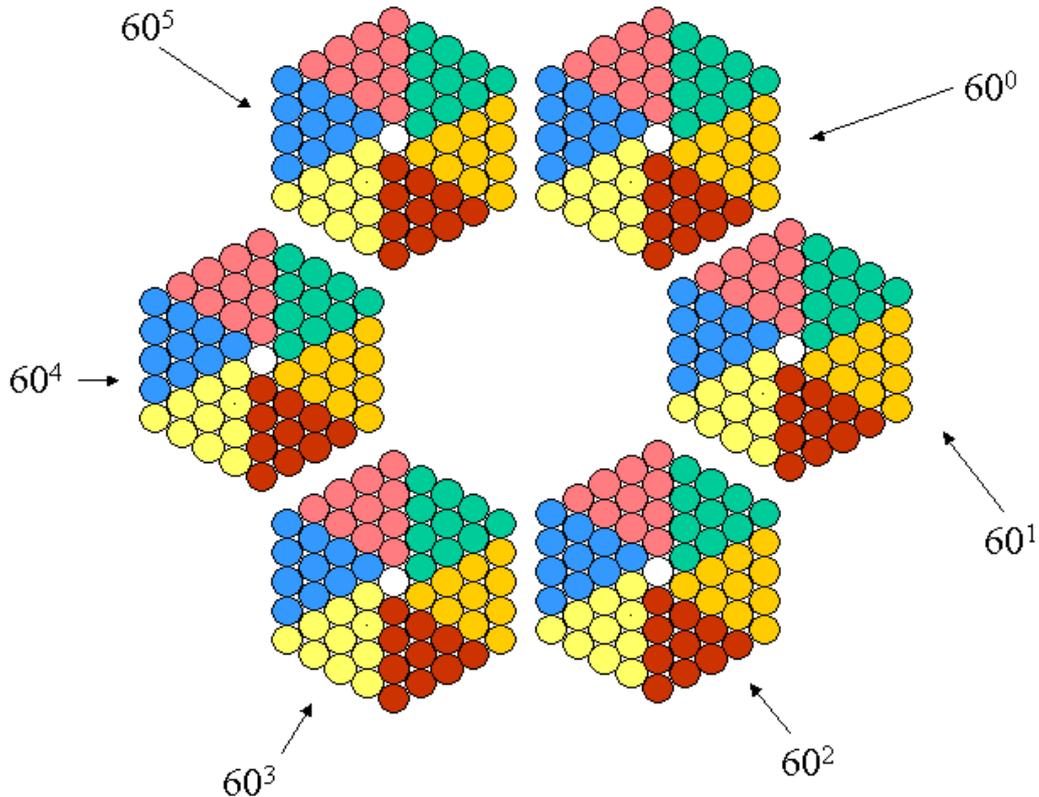

The diagram above shows the weight that would have pieces like stones or tokens placed on the 6 "hexagons". And here is a representation of the number 6''' 27 '' 49 ' (in base 60):

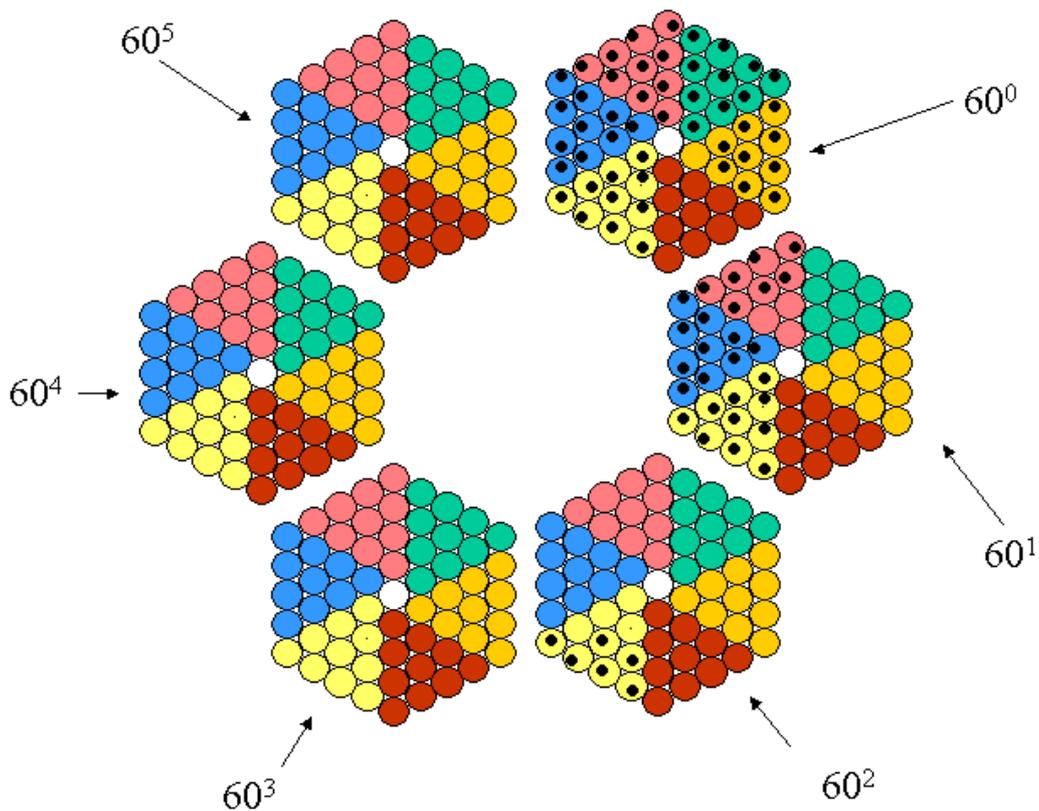

In decimal base this number will be equivalent to:

$6 \times 60^2 + 27 \times 60^1 + 49 \times 60^0 = 23269$

It is clear that with this table one can carry out sums by reducing the pieces of an hexagon filled with a piece of the following hexagon like one makes with the abacuses of decimal base.

**5. - DIVISION OF THE CIRCLE IN 6 AND 12 EQUAL PARTS:**

One draws initially an equilateral triangle.

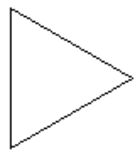

One juxtaposes then 5 other equilateral triangles of the same size as one can see below:

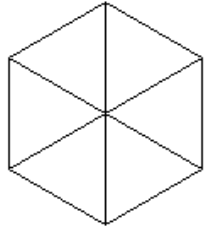

In that way one obtains an hexagon. Then, around the centre of this hexagon one will draw a circle with a radius of half the length of the sides of the equilateral triangles.

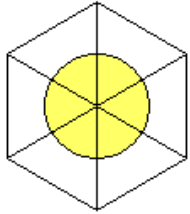

It is clear now that the circle is divided into 6 equal parts, each part being delimited by one of the radiuses of the hexagon. This determines also 6 equal angles whose vertex is the centre of the circle.

If one draws other circles, with the same radius than the first one, around the vertexes of the hexagon, one will obtain the following figure, what one has seen previously, in point 3:

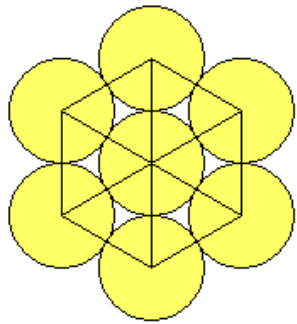

Then one can add equilateral triangles of the same size than the first ones in the following way:

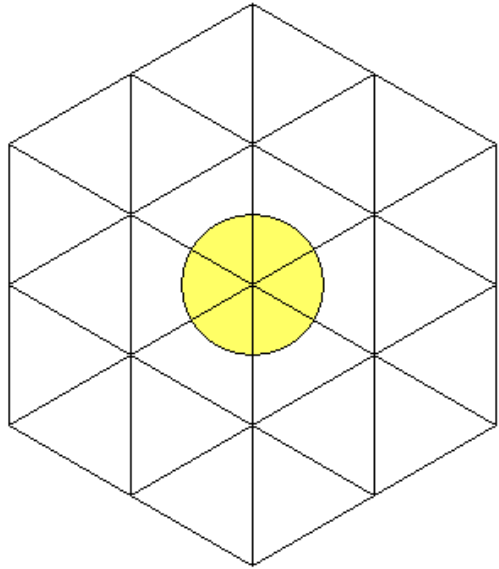

As before, if one draws circles around the vertexes of the equilateral triangles as it is showed below, one obtains again the figure seen in point 3 of this article:

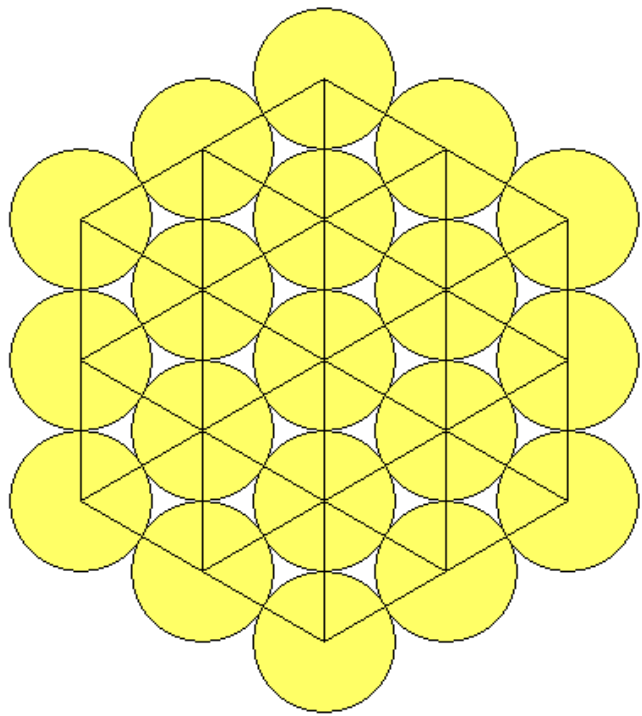

This shows that this construction with the equilateral triangles is linked to the operation that consists in adding orbits of discs as we have seen above.

Now, if one traces three brown radiuses as indicated below, one can observe that the angle of the circle that has been found before will be divided into two equal parts in its turn. This is explained by using symmetries and properties of equilateral triangles.

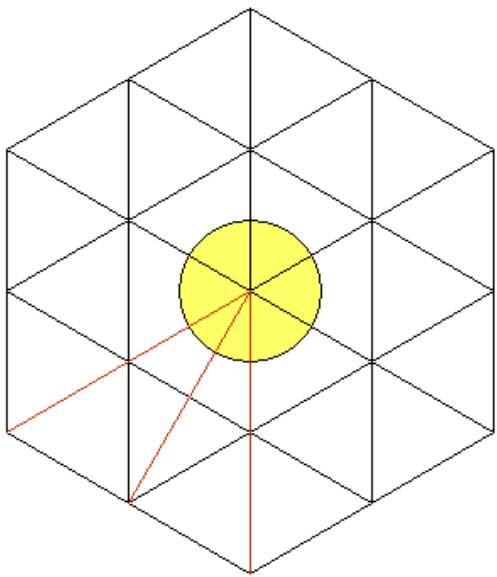

As this process can be applied to the rest of the great hexagon, one can say that the circle of the centre has been divided into 12 equal parts and that 12 equal angles have been found:

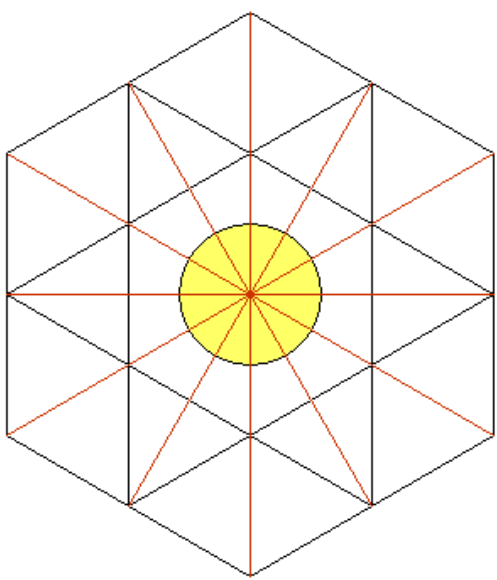

**6. - DIVISION OF THE CIRCLE IN 360 UNEQUAL PARTS BUT RELATIVELY CLOSE TO EACH OTHER**

One can obviously add more equilateral triangles, as it has been done until now, around those that have already been drawn:

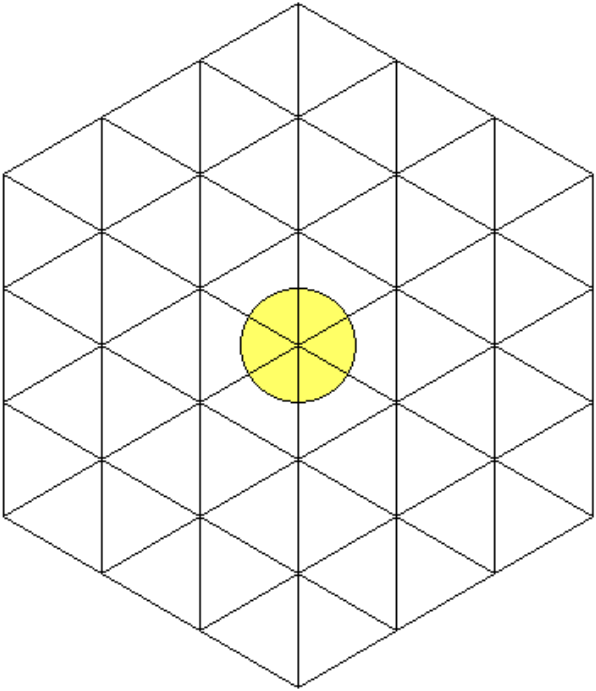

And once again one can divide the 6 initial angles of the circle into other parts, unequal this time, as it can be noticed below:

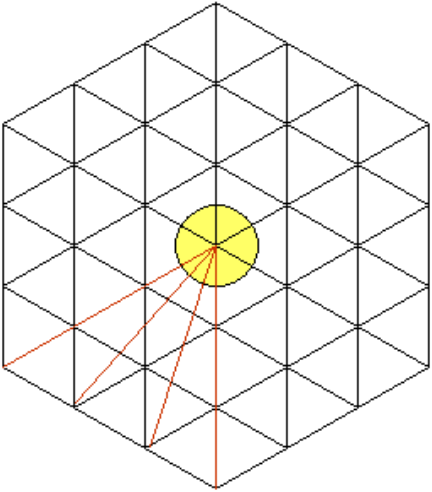

This time one can observe that the initial angle has been divided into 3 parts. As at the beginning there were 6 angles, one can say that the circle has been divided into 3x6=18 parts and that also 18 angles have been found, but unequal ones.

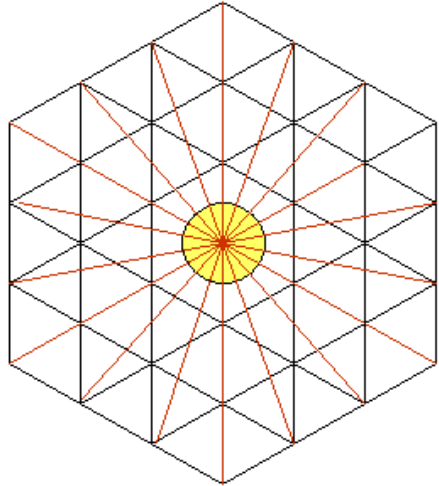

One can continue to add equilateral triangles around the last hexagon as it has been done until now, but one already has been able to notice that there were always 6 large equilateral triangles composed of under-triangles, equilateral ones also, around the centre of the circle. Here there are 6 large triangles distinguished by green lines:

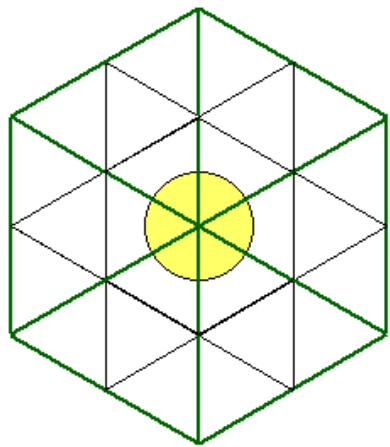

Here are the large triangles with more under-triangles around the circle:

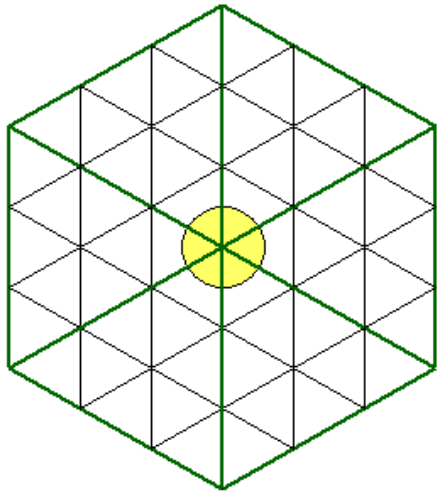

The fact that the equilateral under-triangles, laid out as proposed, find themselves inside other equilateral triangles is easily explained thanks, once again, to symmetries and the properties of the equilateral triangles.

One can thus concentrate on one of the large equilateral triangles in order to see what occurs thereafter:

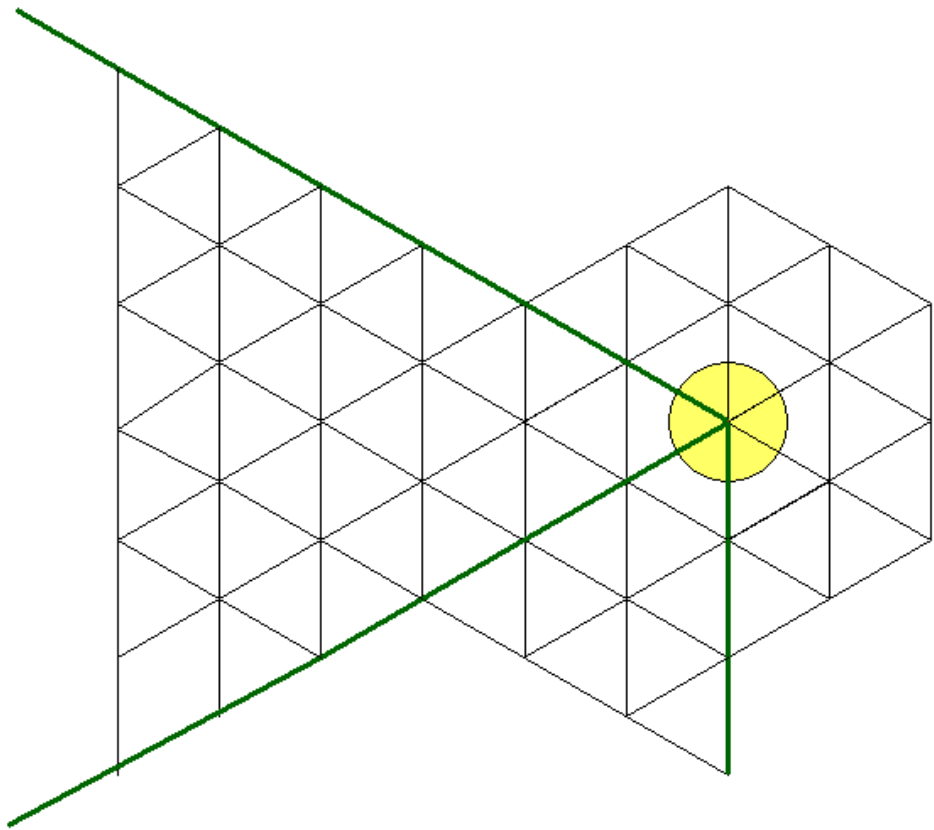

As it can be noticed above, with each new addition of equilateral under-triangles, according to the suggested method, the 6 large triangles, delimited by the green lines, increase by an odd number of triangles. And also, at each stage, one can divide the initial angle into 2 parts, then into 3, then into 4, then into 5, then into 6...:

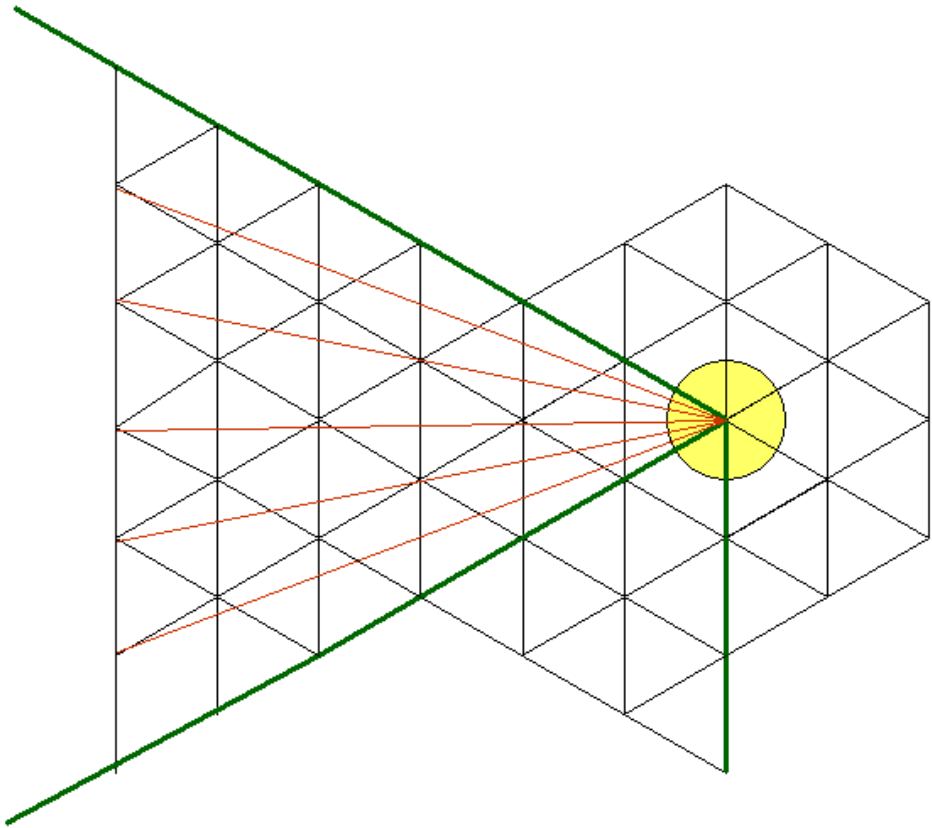

It is shown clearly above that the initial angle has been divided into 6 parts. As there were at the beginning six equal angles, one can say that, in this way, 6x6=36 unequal angles have been found. Moreover, one can conclude that the number of parts which one obtains is a multiple of 6.

If one continues the process, it is clear from now on that one will be able to divide the initial angle into 60 parts. As there were 6 equal angles at the beginning, the circle will have been divided into 60x6=360 unequal parts. In this way we have succeeded in dividing the circle in 360 parts, in a purely geometrical manner. And it is clear from now on that even if one doesn't obtain equal parts, the division of the circle in 360 parts is however linked to the basic properties of the circle and the equilateral triangle. From this idea one can try to obtain a division in 360 equal parts as the ancient Greeks had tried to do.

However, it should be noticed that even if the 360 angles, obtained previously, are unequal, they are nevertheless relatively close to the values of the equal angles that correspond each one to the $360^{th}$ part of the circle. Indeed, one obtains the aforementioned unequal angles by equal divisions of the 6 chords corresponding to the 6 angles obtained at the first division of the circle. In fact, the value of each of these 360 unequal angles, obtained with this method, varies between about 0.833883984 degrees and about 1.10252169 degrees. To observe this better, let us note the angles obtained by the following described manner:

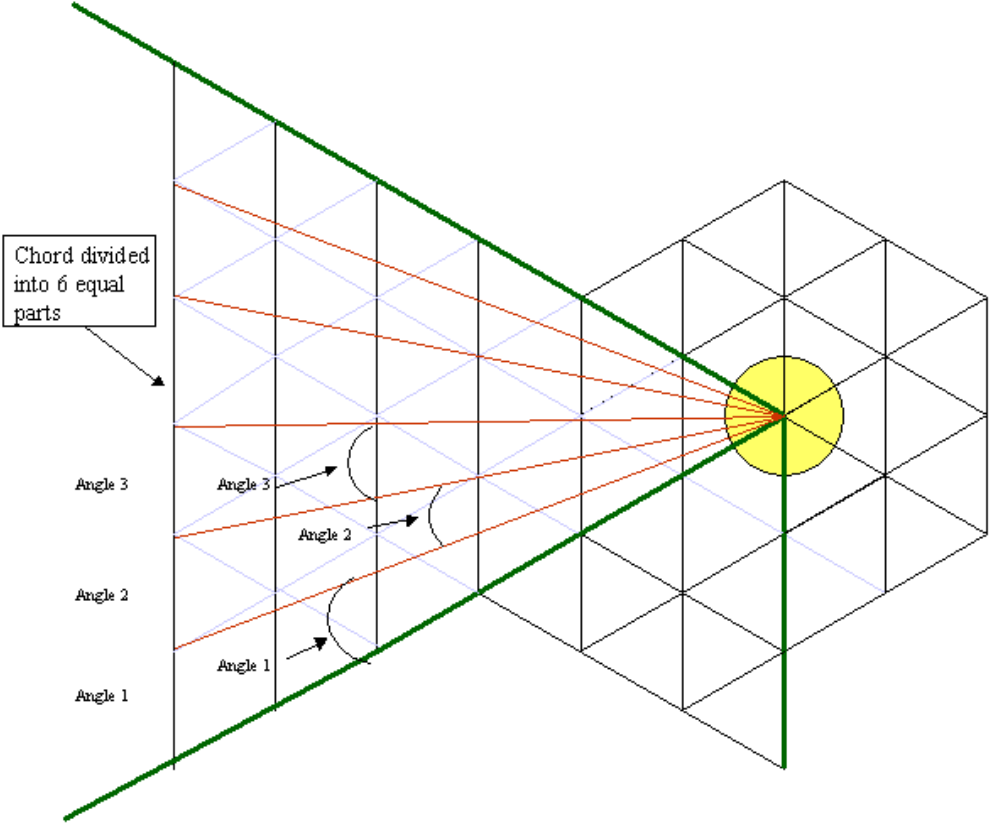

    The diagram that can be seen above shows what the notations "Angle 1", "Angle 2", etc., that are on the left drawing, refer to. These notations refer to the angles indicated by the arrows and which are thus obtained by connecting the centre of the yellow disc to the various points of the chord that has been divided into 6 equal parts.
    Like it has been understood previously, dividing angles with the suggested method returns to divide a chord corresponding to the 6th part of a circle into n equal parts and then connecting to the centre of the central circle each points that delimit the parts of the divided chord. However, in order to divide the circle into 360 parts, one must divide a chord covering the 6th part into 60 parts. But one can also initially work with half of the equilateral triangle and divide the half-chord into 30 parts. In this way one will obtain the following angles:

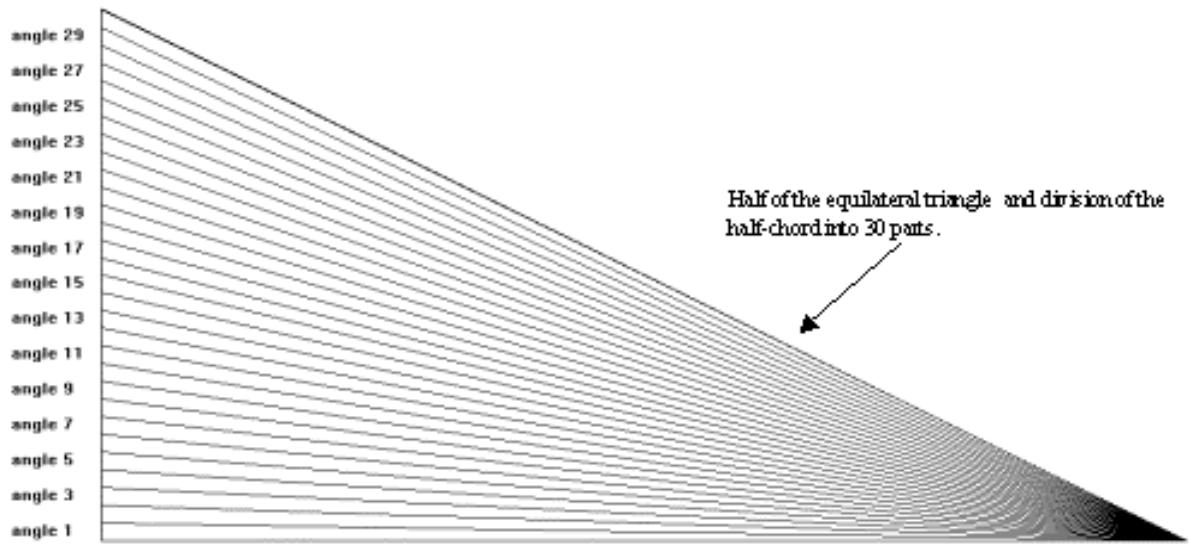

Half of the equilateral triangle and division of the half-chord into 30 parts.

And now we can see the values of these angles which are numbered. On the table below the absolute error corresponds to the difference between the angle obtained and the theoretical angle of 1 degree. The covered angle corresponds to the angle obtained by superposing 60 times each small angle:

| Angle | Value in degrees | Absolute error | Relative error | Covered arc |
|---|---|---|---|---|
| 1 | 1.10252169 | 0.10252169 | 9.299% | 66.1513014 |
| 2 | 1.101705813 | 0.101705813 | 9.232% | 66.1023488 |
| 3 | 1.100077676 | 0.100077676 | 9.097% | 66.0046606 |
| 4 | 1.097644472 | 0.097644472 | 8.896% | 65.8586683 |
| 5 | 1.094416892 | 0.094416892 | 8.627% | 65.6650135 |
| 6 | 1.090409009 | 0.090409009 | 8.291% | 65.4245405 |
| 7 | 1.085638123 | 0.085638123 | 7.888% | 65.1382874 |
| 8 | 1.08012458 | 0.08012458 | 7.418% | 64.8074748 |
| 9 | 1.073891559 | 0.073891559 | 6.881% | 64.4334935 |
| 10 | 1.066964833 | 0.066964833 | 6.276% | 64.01789 |
| 11 | 1.059372514 | 0.059372514 | 5.604% | 63.5623509 |
| 12 | 1.051144779 | 0.051144779 | 4.866% | 63.0686868 |
| 13 | 1.042313587 | 0.042313587 | 4.060% | 62.5388152 |
| 14 | 1.032912388 | 0.032912388 | 3.186% | 61.9747433 |
| 15 | 1.022975834 | 0.022975834 | 2.246% | 61.3785501 |
| 16 | 1.01253949 | 0.01253949 | 1.238% | 60.7523694 |
| 17 | 1.001639548 | 0.001639548 | 0.164% | 60.0983729 |
| 18 | 0.990312561 | 0.009687439 | 0.978% | 59.4187537 |
| 19 | 0.978595178 | 0.021404822 | 2.187% | 58.7157107 |
| 20 | 0.966523906 | 0.033476094 | 3.464% | 57.9914344 |
| 21 | 0.95413488 | 0.04586512 | 4.807% | 57.2480928 |
| 22 | 0.941463658 | 0.058536342 | 6.218% | 56.4878195 |

| 23 | 0.928545037 | 0.071454963 | 7.695% | 55.7127022 |
| 24 | 0.915412887 | 0.084587113 | 9.240% | 54.9247732 |
| 25 | 0.902100007 | 0.097899993 | 10.852% | 54.1260004 |
| 26 | 0.888638004 | 0.111361996 | 12.532% | 53.3182802 |
| 27 | 0.875057188 | 0.124942812 | 14.278% | 52.5034313 |
| 28 | 0.861386496 | 0.138613504 | 16.092% | 51.6831897 |
| 29 | 0.847653424 | 0.152346576 | 17.973% | 50.8592054 |
| 30 | 0.833883984 | 0.166116016 | 19.921% | 50.0330391 |

Obviously, the closer we get to 1 degree, the closer to 60 degrees this covered arc will be. This superposing can help to find the nearest angle to 1 degree. And we observe thus that the 16$^{th}$ and 17$^{th}$ angles are very close to 1 degree.

The method presented in this paper enables to divide the circle into 6 and 12 equal parts. By simple division in 2 with the compass one can divide again the circle into 24, 48, 96, 192 and 384 equal parts.

In order to attain a division in 360 equal parts, at first it would be necessary to divide the circle into 24 parts and then to be able to divide the angles into 3 and 5 parts since 360=24*3*5. Therefore that supposes a method to divide into 3 or 5 the angle corresponding to the 24$^{th}$ part of the circle, i.e. 15 degrees. If one stops at 12 angles at the beginning, it would be necessary first to divide an angle of 30 degrees into 3 parts and then into 5 parts. However, one cannot, using the ruler and the compass, divide most of the time an angle into 3 and 5 parts.

**7.- Star of David :**

The way of surrounding a circle by other circles, as it has been seen at the beginning, makes it possible to draw the star of David:

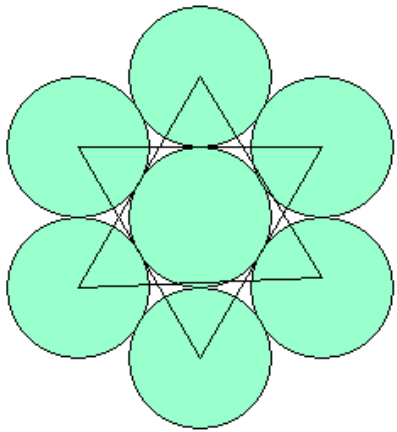

In order to get this it is enough to join suitably the centres of the 6 circles of the first orbit.

**8. - Infinitesimal properties:**

One has seen in point 3 that the more orbits one adds, the closer one gets to an hexagon formed of discs with the same radius. This can be noticed more clearly thanks to the diagram below:

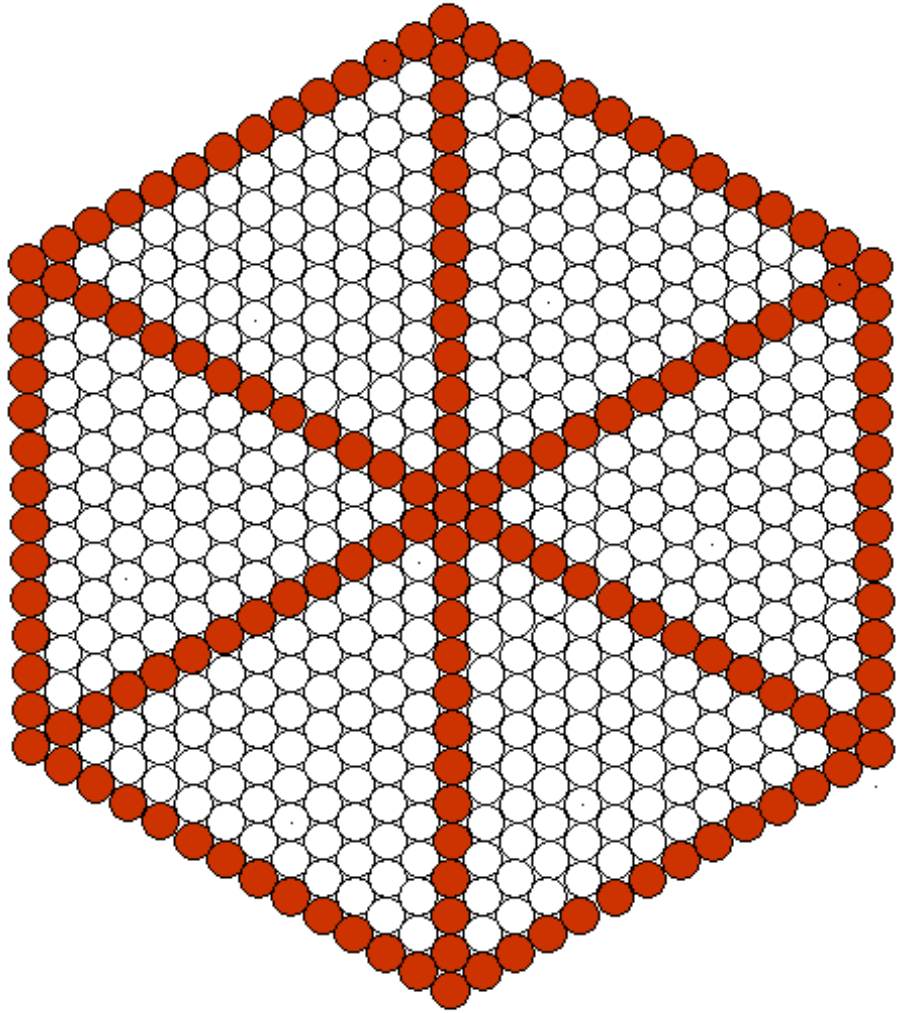

It means that if one continues to add orbits indefinitely while reducing the size of the discs so that the hexagon formed by the discs is inscribed in a circle of radius 1, one will obtain this "limit" form:

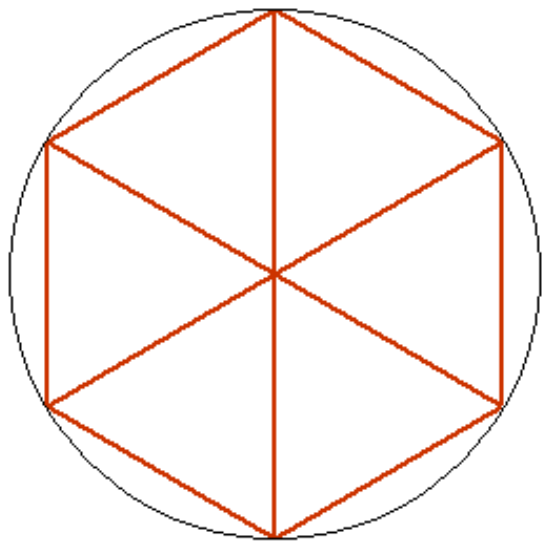

In other words, one obtains an hexagon starting from an infinity of discs laid out in the way that is indicated. And this limit-hexagon will be itself inscribed in a circle...

## 9.- Conclusions

As it has been seen, these geometrical facts show that there is a geometrical link between the circle and the sexagesimal system and that the division of the circle in 360 equal parts, which corresponds to 360 equal angles or degrees, isn't a choice that is completely independent of the geometry because this division can be suggested by the basic properties of the circles and the equilateral triangles. Moreover, the exposed geometrical facts make it possible to represent astronomical data in a way that relates to the movement of the planets around a centre. And one also has been able to represent numbers in base 60 thanks to these properties of the circle and geometrically find the principal divisors of number 60.

It may be that the way of representing time, the division that has been tried in 360 more or less equal parts and the representation of the numbers in base 60 that have been set out in this article played a part during the adoption of the sexagesimal system to represent the numbers, divide the circle and represent time.

## 10. - Comments

This article is the third version of the presentation of these geometrical facts. At present the author is preparing another more complete version with calculations and precisions requested by the readers of this version.

**Version 3 of the text entitled « UNE METHODE NATURELLE POUR DIVISER LE CERCLE EN 360 PARTIES EGALES, CE QUI CORRESPOND A DES DEGRES », written in January 2005.**

**Original title : « UN LIEN GEOMETRIQUE ENTRE LE CERCLE ET LE SYSTEME SEXAGESIMAL »**

**Author: Jaime Vladimir TORRES-HEREDIA JULCA**

**July 2005**

## Complementary notes to the text entitled « A geometrical link between the circle and sexagesimal system » written by Jaime Vladimir TORRES-HEREDIA JULCA :

1.- Using orbits seen in point 3 one can build "hexagons" formed of circles of the same radius in such a way that with six hexagons one will have 360 coloured discs (one doesn't count the circles of the middle). With that one can thus build a calendar as it is seen below:

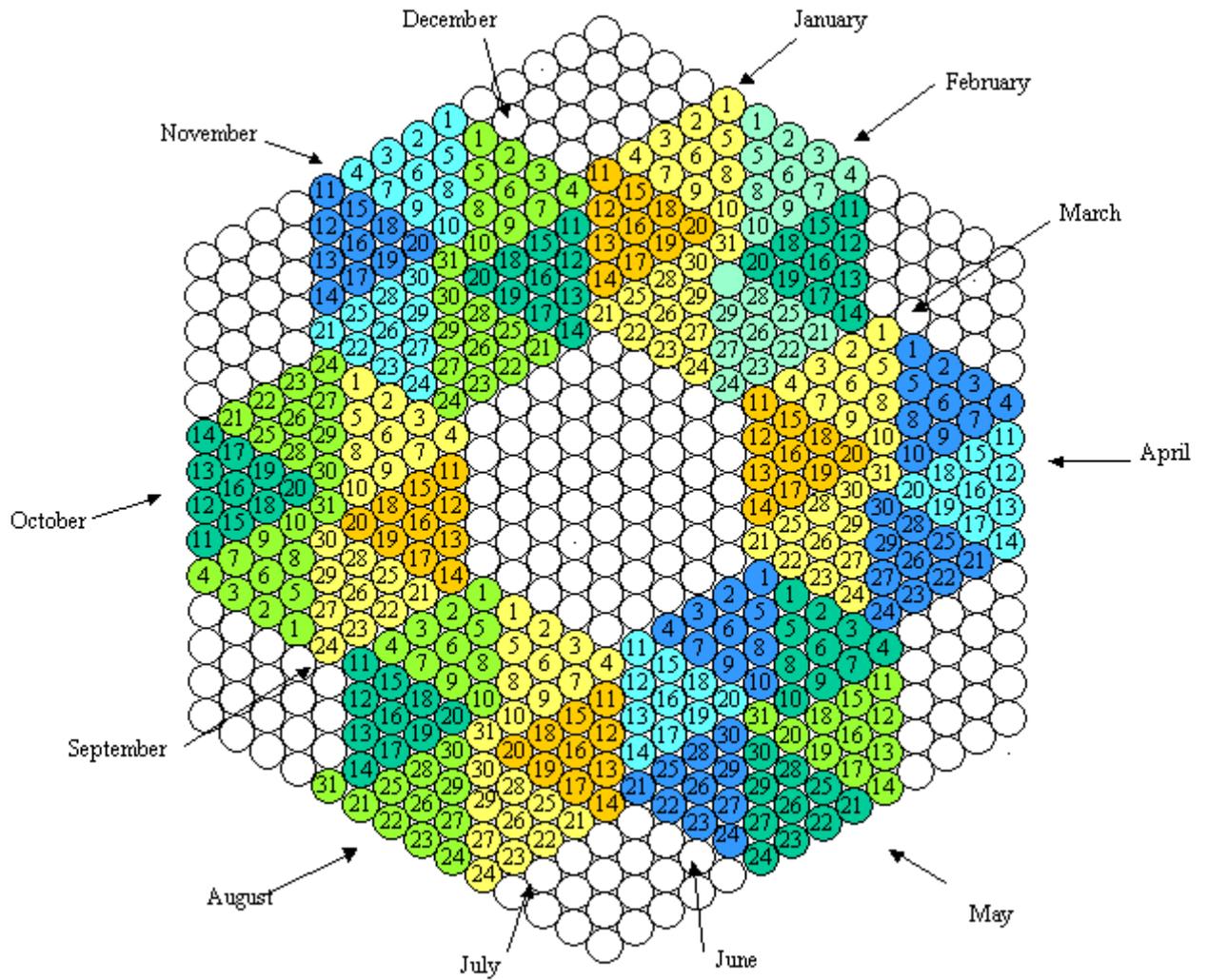

The months are distinguished by the colours yellow, green and blue so that one distinguishes the groups of 10 days for each month.